\documentclass[11pt]{article}
\usepackage[english]{babel}
\usepackage[utf8]{inputenc}
\usepackage{todonotes}
\usepackage{amsmath,amsfonts,amssymb,amsthm,mathrsfs}
\usepackage{amsbsy,mathtools,stmaryrd,upgreek,mathabx}
\usepackage{lmodern}
\usepackage{subcaption}
\renewcommand{\leq}{\leqslant}
\renewcommand{\geq}{\geqslant}

\usepackage[mono=false]{libertine}

\usepackage[backend=biber,
style=alphabetic,]{biblatex}
\usepackage{csquotes}
\newcommand{\la}{\left\langle}
\newcommand{\ra}{\right\rangle}
\newcommand{\norm}[1]{\lVert #1 \rVert}
\newcommand{\abs}[1]{\left| #1 \right|}

\newcommand{\defeq}{\mathrel{\mathop:}=}

\def\R{\mathbb{R}}

\def\N{\mathbb{N}}
\def\Z{\mathbb{Z}}
\def\P{\mathbb{P}}
\def\E{\mathbb{E}}
\def\eps{\varepsilon}
\def\ind{\mathbf{1}}

\newcommand\mydots{\hbox to 1em{.\hss.\hss.}}

\newtheorem{theorem}{Theorem}[section]
\newtheorem{lemma}[theorem]{Lemma}

\theoremstyle{definition}
\newtheorem{definition}[theorem]{Definition}

\theoremstyle{remark}
\newtheorem{remark}[theorem]{Remark}

\usepackage[a4paper,vmargin={3.5cm,3.5cm},hmargin={3cm,3cm}]{geometry}

\usepackage[english]{babel}
\usepackage[colorlinks=true]{hyperref}

\SetSymbolFont{stmry}{bold}{U}{stmry}{m}{n}

\title{Genealogies under logistic growth}
\author{Ruairi Garrett\thanks{Department of Statistics, University of Oxford. Email: \texttt{garrett@stats.ox.ac.uk}} \space and Julio Ernesto Nava Trejo\thanks{Department of Statistics, University of Oxford. Institut für Mathematik, Humboldt Universität zu Berlin. Email: \texttt{julio.ernesto.nava.trejo@hu-berlin.de}}}
\date{\today}

\begin{document}

\maketitle
\begin{abstract}
    We derive the asymptotic behaviour of the genealogy of a logistic branching process in the setting where the equilibrium population size is large. In three regimes for the tail of the offspring distribution we recover the Kingman, \(\text{Beta}(2-\alpha, \alpha)\) and Bolthausen-Sznitman coalescents as a scaling parameter governing the population size is taken to infinity; the deduction goes via the convergence in distribution of a modified lookdown construction. This resolves a question left open in \cite{F25} who studied the same population process forwards in time, showing convergence of the type frequency process to the corresponding \(\Lambda\)-Fleming-Viot process in each regime. Further, when the type space is the continuous torus and the individuals mutate their location as independent Brownian motions we recover the Brownian Spatial \(\Lambda\)-coalescent. 
\end{abstract}
\section{Introduction}

A primary objective of mathematical population genetics is to identify the genealogical relationship between individuals sampled from a population that evolves according to some prescribed dynamics. This becomes challenging in the typical situation when the population size is not reversible and backwards in time descriptions of the genealogy are not available. Here we use a lookdown construction to overcome this issue and identify the genealogy of the logistic branching processes with large carrying capacity, which has been studied as a model of a competing population. This extends the work of \cite{F25} so we begin by briefly presenting their main results. 

\subsection{Logistic branching and \((\Lambda, \mathcal{Q})\)-Fleming-Viot processes}
Fix a scaling parameter \(K\). Consider a population of individuals each carrying a neutral genetic type in some compact Polish space \(E\). This type mutates in \(E\) according to a c\`adl\`ag Markov process with generator \(Q^K\) independently of everything else. Individuals give birth at rate \(b\) to a random number of children with law \((p_\ell)_{\ell\geq 1}\), each of which inherits the genetic type of its parent the instant before birth. Letting \(N^K(t)\) denote the size of the population at time \(t\geq 0\), the rate at which each individual dies is given by the competition term \(d + \frac{c}{K}N^K(t)\), where \(c,d>0\) are constants. We record the state of the population by its empirical measure: let \((Y^K_i(t))_{i=1}^{N^K(t)}\) denote the types of individuals alive at time \(t\) in this population, and set \[\nu^K(t) = \sum_{i = 1}^{N^K(t)}\delta_{Y^K_i(t)} \] which belongs to \(\mathcal{M}_{F}(E)\), the space of finite measures on \(E\). Since it will be useful later we characterize \(\nu^K\) as the solution of a martingale problem. Let \(g:E\to \R\) be such that \(0\leq g\leq 1\) and \(g\in \mathcal{D}(\mathcal{Q}_K)\), and set \(F_{g}(\nu) = \exp(\langle \log g, \nu \rangle) = \prod\limits_{x\in \nu}g(x)\), under the common abuse of notation in which counting measures are treated as sets. Let \(\mathcal{D}\) be the linear span of the functions \(F_{g}\). 

\begin{definition}[Logistic branching process with neutral mutation]\label{DEF:logisticbranchingprocess}
    The logistic branching process \(\left(\nu^K(t)\right)_{t\geq 0}\) is the unique Markov process taking values in \(\mathcal{M}_F(E)\) which solves the martingale problem \((\mathcal{L}^K, \mathcal{D})\), where for \(\mathcal{L}_K\) acts on functions of the form \(F_g\) by
    \begin{align}
        \label{eq:gen_K_BGW}
    \notag \mathcal{L}_{K}F_{g}(\nu)\coloneqq &  \sum_{\ell=1}^{\infty}\int_{E}bp_\ell\left(\exp(\la \log g, \nu\ra +\ell\log g(x))-\exp(\la \log g, \nu\ra)\right) \nu(dx)\\
    & + \left(d+\frac{c}{K}\la 1, \nu\ra\right)\int_{E}\left(\exp(\la \log g, \nu\ra -\log g(x))-\exp(\la \log g, \nu\ra)\right) \nu(dx). 
    \\& + \sum_{x\in\nu} \mathcal{Q}^Kg(x) \prod_{y\neq x} g(y)
    \end{align}
\end{definition}

At this level of generality the above martingale problem may not be well posed, but it is for all conditions on \((p_\ell)_{\ell\geq 1}\) appearing in our results. Whenever \(b\sum_{\ell}\ell p_\ell > d\), the population size behaves like a supercritical branching process when small but becomes subcritical once the number of individuals exceeds order \(K\) (at least when the offspring distribution has finite mean). Consequently, the population size stabilizes on order \(O(K)\) and we focus on the large \(K\) regime. The mutation mechanism is included purely to formulate Theorem \ref{thm:spatial} below; for the purpose of our main result Theorem \ref{thm:main}, it can be ignored.

To state the main result of \cite{F25} we recall the definition of the \(\Lambda\)-Fleming-Viot processes which appeared first implicitly in \cite{DK99} and were formally introduced by Bertoin and Le Gall \cite{Bertoin2003} as a generalization of the classical Fleming–Viot process. Further, for the generator \(\mathcal{Q}\) of a c\`adl\`ag Feller Markov process we also recall the definition of the \((\Lambda, \mathcal{Q})\)-Fleming-Viot  processes (see e.g. \cite{Xi09}). These are Markov processes taking values in the space \(\mathcal{M}_1(E)\) of probability measures on \(E\), also characterized as the unique solution to a martingale problem. Let \(\Lambda\) be a finite measure on \([0, 1]\), and decompose \(\Lambda = a\delta_0 + \Lambda_0\), where \(\Lambda_0\) has no atom at \(0\). For \(f:E^n \to \R\) bounded and measurable, consider test functions of the form \[G_f(\rho) = \int_{E^n}f(\mathbf{x})\rho^{\otimes n}(d\mathbf{x}).\]
\begin{definition} The \((\Lambda, \mathcal{Q})\)-Fleming-Viot process is the Markov process with values in \(\mathcal{M}_F(E)\) and generator \(\mathcal{L}_{\Lambda, \mathcal{Q}}\) acting on test functions of the form \(G_f\) via \[\mathcal{L}_{\Lambda, \mathcal{Q}} G_f(\rho) = a \mathcal{L}_0 G_f(\rho) + \mathcal{L}_{\Lambda_0}G_f(\rho) + \mathcal{L}_\mathcal{Q} G_f(\rho), \] where \begin{align*}
    \mathcal{L}_0 G_f(\rho) &= \sum_{1 \leq i < j \leq n} \int_{E^n}\Big(f(x_1, \mydots, x_i, \mydots, x_i, \mydots, x_n) - f(x_1, \mydots, x_i, \mydots, x_j, \mydots, x_n)\Big)\rho^{\otimes n}(d\mathbf
    {x}),\\
    \mathcal{L}_{\Lambda_0}G_f(\rho) &= \int_{[0,1]}\int_E G_f\left((1-u)\rho + u\delta_x\right)\rho(dx)\frac{\Lambda(du)}{u^2},\\
    \mathcal{L}_\mathcal{Q} G_f(\rho) &= \sum_{i = 1}^n\int_{E^n}\mathcal{Q}_if(\mathbf{x})\rho^{\otimes n}(d\mathbf{x}),
\end{align*}
    and \(\mathcal{Q}_i\) denotes the generator \(\mathcal{Q}\) acting on the \(i\)-th coordinate. The case \(\mathcal{Q} = 0\) is the \(\Lambda\)-Flemin-Viot process, and in the special case \(\Lambda = \frac{1}{N_e}\delta_0\) this is the Fleming-Viot diffusion with effective population size \(N_e\) and mutation mechanism \(\mathcal{Q}\). 
\end{definition}


We return to the logistic branching process. Letting \(x_0 \in E\) be arbitrary, define \[\rho^K(t) \coloneqq \begin{cases}
		\frac{\nu^K(t)}{\langle 1, \nu^K(t)\rangle} & \text{ if}\quad \langle1, \nu^K(t)\rangle > 0, \\
		\delta_{x_0} & \text{ otherwise.} 
	\end{cases}\] 
    For a space \(S\), recall that \(D_S[0, T]\) denotes the space of c\`adl\`ag paths \([0, T]\to S\) equipped with the Skorokhod \(J_1\) topology. Whenever \(S\) is a finite set it is given the discrete metric without further comment. The main result of \cite{F25} is to establish that the measure-valued process \(\rho^K\) of the types in \(E\) converges as \(K \to \infty\) to a corresponding \(\Lambda\)-Fleming-Viot process, where the parameter \(\Lambda\) is determined by the offspring law \((p_\ell)_{\ell\geq1}\). 
\begin{theorem}(\cite{F25}, Theorems 2, 3, 4)\label{Thm:raphael} Specialise to the case \(\mathcal{Q}_K = 0\) (no mutation).
     \begin{enumerate}
       \item \textbf{Finite variance regime.} Suppose for some \(\delta > 0\) that \(\sum_{\ell\geq1}p_{\ell} \ell^{2+\delta} < \infty\). Set \(m = \sum_{\ell\geq 1} \ell p_\ell\) and \(m^{(2)} = \sum_{\ell}\ell^2p_\ell\). Assume that \(bm-d > 0\) and set \(n_* = \frac{bm-d}{c}\). Set \(n^K(t) = \frac{1}{K}N^K(Kt)\). If \(n^K(0) \to n_*\) and \(\rho^K(0) \to \rho(0)\) in probability, then \[\Bigl(\rho^K(Kt), t\in[0, T]\Bigl) \Rightarrow \Bigl(\rho(t), t\in[0, T]\Bigl)\] in distribution in \(D_{\mathcal{M}_1(E)}[0, T]\), where \(\rho\) is the Fleming-Viot diffusion with effective population size \[N_e = \frac{n_*}{b(m + m^{(2)})}.\]
       \item \textbf{The \(\boldsymbol{\alpha}\)-stable regime.} Suppose \(p_\ell\sim \frac{p_0}{\ell^{1+\alpha}}\) as \(\ell \to \infty\) for \(\alpha \in (1, 2)\). Set \(n^K(t) = \frac{1}{K}N^K(K^{\alpha-1}t)\). If \((n^K(0), K \geq 1)\) is tight in \((0, \infty)\) and \(\rho^K(0)\to \rho(0)\) in probability, then \[\Bigl(\rho^K(K^{\alpha-1}t), t\in[0, T]\Bigl) \Rightarrow \Bigl(\rho(t), t\in[0, T]\Bigl)\] in distribution in \(D_{\mathcal{M}_1(E)}[0, T]\), where \(\rho\) is the \(\Lambda\)-Fleming-Viot process with \[\Lambda(du) = \frac{bp_0}{n_*^{\alpha-1}}(1-u)^{\alpha-1}u^{1-\alpha}du.\]
       \item \textbf{Neveu's branching process with logistic competition.} Suppose \(p_\ell\ell^2 \to p_0\) as \(\ell \to \infty\). In this case redefine \(n^K(t) = \frac{1}{K\log(K)}N^K(t)\). If \((n^K(0), K \geq 1)\) is tight in \((0, \infty)\) and \(\rho^K(0)\to \rho(0)\) in probability, then \[\Bigl(\rho^K(t), t\in[0, T]\Bigl) \Rightarrow \Bigl(\rho(t), t\in[0, T]\Bigl)\] in distribution in \(D_{\mathcal{M}_1(E)}[0, T]\), where \(\rho\) is the \(\Lambda\)-Fleming-Viot diffusion with \[\Lambda(du) = bp_0du.\] 
   \end{enumerate}
\end{theorem}

 It is also shown in \cite{F25} that the rescaled population size stays near a certain carrying capacity, we will recall precise statements when necessary and make heavy use of this shortly. 
\subsection{Coalescents, spatial coalescents and Brownian spatial \(\Lambda\)-coalescents}
The \((\Lambda,\mathcal{Q})\)-Fleming-Viot processes have a well known genealogical structure, which we now recall. A \(k\)-coalescent is a c\`adl\`ag Markov process taking values in the set \(\mathcal{P}_k\) of partitions of \([k] \coloneqq \{1,\ldots,k\}\) in which the only transitions consist of merging blocks. It is label invariant if its law is unchanged under relabelling the initial blocks. Label invariant coalescents with no simultaneous mergers are characterized by a finite measure \(\Lambda\) on \([0, 1]\), leading to the class of \(\Lambda\)-coalescents introduced independently by Donnelly and Kurtz, Pitman, and Sagitov \cite{DK99,pitman1999,sagitov1999} as a generalization of the famous Kingman's coalescent \cite{kingman1982genealogy}. For our purpose these are Markov processes taking values in \(\mathcal{P}_k\), the set of partitions of \([k]\), with transition rates defined as follows: for \(\pi, \pi' \in \mathcal{P}_k\) such that \(|\pi| = n\) and \(\pi'\) is obtained from \(\pi\) by merging exactly \(j\) blocks, the transition rate is
\[
\lambda_{\pi,\pi'} = \int_0^1 u^{j-2}(1 - u)^{n - j} \Lambda(du),
\]
and \(\lambda_{\pi,\pi'} = 0\) otherwise.
\begin{remark}
    \(\Lambda\)-coalescents are usually defined to take values in partitions of \(\N\). Here we are only concerned with the \(k\)-coalescents relating samples of \(k\) individuals from our population model, hence the restriction to \([k]\).
\end{remark}
It is well known that the genealogy of the \(\Lambda\)-Fleming-Viot process is given by the \(\Lambda\)-coalescent \cite{DK99,Bertoin2003}. Particularly relevent here is the pathwise duality established in \cite{DK99} which will be discussed soon. With this in mind, the result of Theorem \ref{Thm:raphael} strongly suggests that the genealogy of the prelimiting population models should also converge to the corresponding \(\Lambda\)\nobreakdash-coalescent. Specifically, the genealogy of the population in regime 1 should scale to Kingman’s coalescent \cite{kingman1982genealogy}, in regime 2 it should scale to the \(\mathrm{Beta}(2-\alpha, \alpha)\)-coalescents for \(\alpha \in (1,2)\), and in regime 3 to the Bolthausen–Sznitman coalescent \cite{Bolthausen1998}, which is the special case when \(\Lambda\) is a multiple of the Lebesgue measure on \([0, 1]\). The proof of this conjecture stated in \cite{F25} is the primary objective of this article (see Theorem \ref{thm:main}). These coalescents are collectively known as the Beta-coalescents, since the Beta\((2-\alpha, \alpha)\) distribution converges weakly to \(\delta_0\) as \(\alpha\) goes to \(2\), and to the uniform distribution on \([0,1]\) as \(\alpha\) goes to \(1\). 

In the special case where \(E = \R^d/\Z^d\) is the flat torus and the spatial motion (mutation mechanism) converges on the appropriate timescale to Brownian motion\footnote{See \cite{K25} Section 1.7.2 for discussion on generalities.} a much stronger result can be obtained at almost no additional expense; namely convergence of the genealogy (including the spatial locations of the ancestral lines) to the so-called Brownian spatial coalescent recently introduced in \cite{K25} primarily as a universal object in its own right, but additionally as the genealogy of the \((\Xi, \Delta)\)-Fleming-Viot process on the torus with Brownian mutation. We recall a few informal definitions. Let \(\mathcal{P}\) denote the set of partitions of finite subsets of \(\N\).

\begin{definition}[Spatial coalescent, Brownian spatial coalescent] A spatial coalescent on \(E\) is a c\`adl\`ag\footnote{In the topology obtained by identifying \(\mathbf{X}_t\) with a point in \(E^{\abs{\Pi_t}}\).} Feller Markov process \(\left(\Pi_t, \mathbf{X}_t\right)_{t\geq 0}\) with \(\Pi_t\in \mathcal{P}\) denoting the collection of lineages at time \(t\), \(\mathbf{X}_t:\Pi_t\to E\) assigning a spatial location to each lineage, and such that the only transitions in the first coordinate consist of merging partition elements. 
In the case when \(E\) is the \(d\)-dimensional torus, a spatial coalescent is a Brownian spatial coalescent if, conditional on the topology of the genealogical forest and the times and locations of all mergers, lineages follow independent Brownian bridges along the branches. 
\end{definition}

In the case of no simultaneous mergers, Theorem 1.7 of \cite[Theorem 1.7]{K25} states that label invariant Brownian spatial coalescents are uniquely characterised by a collection of measures \(\mathbf{\mu} = (\mu_{n,k}\in \mathcal{M}_F(E):2\leq k \leq n)\), and \cite[Theorem 1.12]{K25} states that they statisfy an appropriate form of sampling consistency if and only if they are of the form \(\nu_{n, k}(dz) = \lambda_{n, k}dz\) with \(\lambda_{n, k}\) the collision rate in a non-spatial \(\Lambda\)-coalescent. The Brownian spatial coalescent corresponding to this family of measures is then called the Brownian spatial \(\Lambda\)-coalescent. For more particulars, including a description of the law of the Brownian spatial \(\Lambda\)-coalescent, see \cite{K25}.

\subsection{Main results}
In our setting, we describe the genealogy of the population by the ancestral partition process \(\Pi^K\) as follows. 
\begin{definition}[Ancestral partition process, spatial ancestral partition process]
    Let \(r_K\) be either \(K, K^{\alpha-1}\) or \(1\) according to the time rescaling of Theorem \ref{Thm:raphael}. Sample \(k\) individuals uniformly at random without replacement from \(\nu^K(r_KT)\). If there are not enough individuals to sample make an arbitrary choice: the probability of this is vanishing as \(K \to \infty\). Define an equivalence relation \(\sim_{t, K}\) on \([k]\) by declaring \(i \sim_{t, K} j\) if samples \(i\) and \(j\) share a common ancestor at time \(r_KT - t\) (interpreting this in the manner that makes \(\Pi_t\) c\`adl\`ag). For each \(t \in [0, r_KT]\), let \(\Pi^{K}(t)\) be the partition of \([k]\) induced by \(\sim_{t, K}\). Further, as appropriate we denote \(\mathbf{X}^K_t:\Pi^K_t \to E\) the map assigning each block the location of its ancestor at time \(r_KT-t\), modified such that \(\mathbf{X}^K\) has c\`adl\`ag paths. 
\end{definition}
Note that \(\Pi^K\) depends on the choice of \(k\), we suppress this to avoid cluttering the notation. Also note that \(\Pi^K(t)\) is defined for \(t \in [0, r_KT]\). The following is our main result, showing in each of the cases of Theorem \ref{Thm:raphael} that the genealogy \((\Pi^K(r_{K}t), t\in[0, T])\) converges to the corresponding \(\Lambda\)-coalescent.

\begin{theorem}\label{thm:main} In each regime here keep the notation and conditions of the corresponding regime of Theorem~\(\ref{Thm:raphael}\). In each case define \(\Pi^K\) as above. 
    \begin{enumerate}
       \item \label{thm:main.1} \textbf{Kingman's coalescent.} The rescaled ancestral partition process \((\Pi^K(Kt), t \in [0, T])\) converges in distribution in \(D_{P_k}[0, T]\) as \(K \to \infty\) to Kingman's coalescent with effective population size \[N_e = \frac{n_*}{b(m + m^{(2)})}.\]
       \item \label{thm:main.2} \textbf{Beta\(\boldsymbol{(2-\alpha, \alpha)}\)-coalescent.} The rescaled ancestral partition process \mbox{\((\Pi^K(K^{\alpha-1}t), t \in [0, T])\)} converges weakly in \(D_{P_k}[0, T]\) as \(K \to \infty\) to a Beta\((2-\alpha, \alpha)\) coalescent, which is to the \(\Lambda\)-coalescent with \[\Lambda(du) = \frac{bp_0}{n_*^{\alpha-1}}(1-u)^{\alpha-1}u^{1-\alpha}du.\]
       \item \label{thm:main.3} \textbf{Bolthausen-Sznitman coalescent}. The ancestral partition process ~\((\Pi^K(t), t \in [0, T])\) converges (with no time rescaling) in distribution in \(D_{P_k}[0, T]\) as \(K \to \infty\) to a Bolthausen-Sznitman coalescent. This is the \(\Lambda\)-coalescent with \[\Lambda(du) = bp_0du.\] 
   \end{enumerate}
\end{theorem}

\begin{theorem}\label{thm:spatial} Suppose \(E\) is the flat torus, and assume that the sequence of processes with generators \(r_K Q_K\) converge in distribution to Brownian motion on compact time intervals. Assume further to the assumptions of Theorem \ref{Thm:raphael} that \(\rho(0)\) is distributed as the stationary distribution of the \((\Lambda, \Delta)\)-Fleming-Viot process. Then in each of the three regimes on \((p_\ell)_{\ell\geq 1}\), the spatial ancestral partition process converges in distribution for Skorokhod's topology to the corresponding Brownian spatial \(\Lambda\)-coalescent.
    
\end{theorem}


The key tools at play here are the so-called lookdown constructions first introduced in \cite{DK96,DK99} which have proven to be a powerful technique for the construction of branching processes in a manner that retains genealogical information, even under passage to various scaling limits. Of particular interest is the technique presented in \cite[Section~5]{DK99} where the \(\Lambda\)-coalescents are derived by tracing the genealogy backward in time in a Moran-type model. This yields a pathwise duality between the \(\Lambda\)-coalescent and the \(\Lambda\)-Fleming–Viot process, analogous to the pathwise duality established in \cite[Section~3]{DK96} for the Fleming–Viot process and Kingman’s coalescent. 

The basic idea of a lookdown-type construction is to assign the particles a `level' (in our case a positive integer) that orders them in a manner affecting how they reproduce. In particular, the lower a particle's level the higher the rate at which it will give birth. Death events will always remove the particle with the highest level, so the levels order the particles by time of death. In this construction each particle represents a `line of descent' rather than an individual in the population model. The surprising fact is that this can be arranged in way that maintains the exchangeability of the types for all fixed times, so long as the initial distribution is exchangeable. The result is that to fully understand the distribution of the genealogy of a finite sample one need only study the genealogy of the particles carrying the lowest levels, about which explicit computations can be made. 

In Section \ref{sec:Lookdown} we will introduce a lookdown construction of the branching process with logistic regulation, which we will denote by \(X^{K}(t) = (X_1(t), \ldots, X_{N^{K}(t)}(t)) \in \bigcup_{n \in \mathbf{N} }E^n\).  In particular, we will show that the law of the empirical measure induced by this vector agrees with that of \(\nu^K\). We will also recall the limiting particle system \(X\) to which we prove the \(X^K\) converge as \(K \to \infty\). In order to formulate the following theorem it is convenient to view \(X^K(t)\) as an element of \(E^\infty\) by extending the vector. We do this more carefully in Section \ref{sec:Conv}. We will deduce our main result from a minor variant of Theorem \ref{thm:fwd} below. This is analogous to the work of \cite{F25} but at the level of the lookdown construction rather than the measure valued process. The proof relies heavily on the techniques of that article.

\begin{theorem}\label{thm:fwd} Keep the notation of Theorem \(\ref{Thm:raphael}\). In each subpart impose the corresponding conditions of Theorem \(\ref{Thm:raphael}\) and let \(X^K\) be the corresponding lookdown construction as introduced in Section \ref{sec:Lookdown}. Suppose also that the \(E\)-valued process with generator \(r_K\mathcal{Q}_K\) converges in law on compact time intervals to the process with generator \(\mathcal{Q}\).
    \begin{enumerate}
       \item \textbf{Fleming-Viot.} \label{Prop:Kingman}
    The processes \(\left(X^K(Kt)\right)_{t\in [0, T]}\) converge weakly in \(D_{E^\infty}[0, T]\) as \(K\to \infty\) to a lookdown representation of the Fleming-Viot process with mutation mechanism \(\mathcal{Q}\) and effective population size \[N_e = \frac{n_*}{b(m+m^{(2)})}.\] 
       \item \textbf{Beta-Fleming-Viot.} \label{Prop:alpha} The processes \(\left(X^K(K^{\alpha-1}t)\right)_{t\in [0, T]}\) converge weakly in \(D_{E^\infty}[0, T]\) as \(K\to \infty\) to a lookdown representation of the \((\Lambda, \mathcal{Q})\)-Fleming-Viot process with \[\Lambda(du) = \frac{bp_0}{n_*^{\alpha-1}}(1-u)^{\alpha-1}u^{1-\alpha}du.\]
       \item \textbf{Dual of the Bolthausen-Sznitman coalescent}. The processes \(\left(X^K(t)\right)_{t\in [0, T]}\) converge weakly in \(D_{E^\infty}[0, T]\) as \(K\to \infty\) as \(K \to \infty\) to a lookdown representation of the \((\Lambda, \mathcal{Q})\)-Fleming-Viot process with \[\Lambda(du) = bp_0du.\]
   \end{enumerate}
\end{theorem}
As we will see the data of the lookdown construction events encodes all of the genealogical information of the population, so convergence of the lookdown construction will suffice to conclude convergence in law of the ancestral partition process.

\subsection{Related work and alternative approaches}

After completing this work, we became aware of significant overlap with a forthcoming paper by Andr\'e, Foutel-Rodier, and Schertzer \cite{AFS25}, who study the genealogy of a similar model via spinal decomposition techniques, building on the work of Bansaye \cite{Bansaye2024, FoutelRodier2023}. In the present setting, the lookdown framework provides a direct full description of the genealogies of the models under study. However, the genealogies of more general structured populations remain an important open problem and appear to present a challenge to the lookdown approach. The spinal decomposition approach may help overcome some of these difficulties and pave the way for a more systematic study of genealogies in general structured populations.

The appearance of Beta-coalescents from models with \(\alpha\)-stable branching now has some history. Of particular relevance to this article is \cite{SCHWEINSBERG2003}, who shows (among other cases) that the Beta-coalescents emerge from a model in which each individual in each generation produces a number of offspring with law \(\xi\), where \(\P(\xi \geq x)\sim \frac{C}{x^\alpha}\), \(\alpha \in (1, 2)\). From this pool of potential offspring some large number \(N\) are selected to form the next generation. In this model, the genealogy approximates the Beta\((2-\alpha, \alpha)\) coalescent as \(N\to \infty\). Here we observe the same effect, but instead of fixing the population size by sampling it is held roughly constant by competition. 

In the direction of \cite{F25} on which the present article is based, Etheridge and March \cite{E&M1991} showed that the Dawson–Watanabe superprocess conditioned to have fixed population size is a Fleming–Viot process. Forien achieves this effect by competition rather than conditioning. The result of \cite{E&M1991} was later extended by Perkins \cite{Perkins1992}. In \cite{Birkner2005} the \(\text{Beta}(2-\alpha, \alpha)\)-coalescents are identified in the non-spatial setting as the genealogy of continuous state branching processes. There, the authors make a time change of the \(\alpha\)-stable branching process using a function of the population size to make the genealogies Markovian, and show that \(\alpha\)-stable branching mechanisms are the only class of continuous state branching processes for which this can be done. In a related approach, \cite{Caballero2024} study the type frequency processes arising from the ratio of two continuous-state branching processes with possibly different reproduction mechanisms, which allows them to construct a homeomorphism between the spaces of continuous-state branching processes and the \(\Lambda\)-coalescents.

More broadly, there appear to be surprisingly few rigorous results identifying genealogies of populations experiencing competition. In \cite{Cheek2022}, Cheek studies a population evolving as a Markov branching process with generalized logistic growth, and manages to recover explicit formulae for the genealogy of a sample taken while the population is still growing superlinearly, working also in the present setting in which the carrying capacity is scaled to infinity. In a similar direction, Dai Pra and Kern \cite{DaiPraKern25} study the genealogy of a two type logistic branching processes with selection, up to the point where the process almost reaches its carrying capacity. Their approach is based on couplings with branching processes, as in \cite{Cheek2022}, together with the use of ancestral selection graphs to recover the genealogy. In their forward in time analysis, they are able to recover the Gillespie–Wright–Fisher diffusion in the limit of infinite carrying capacity by applying the theory of stochastic equations pushed toward a stable manifold by a large drift.

\subsection{Outline of the article}
The proofs of regime  2 and 3 of Theorems \ref{Thm:raphael} and \ref{thm:fwd} rely on an adaptation of the stochastic averaging principle of Kurtz \cite{avg}. We recall this adaptation from \cite{F25} in Appendix~\ref{averagingappendix}. The proof of regime 1 is simpler and follows from standard generator calculations and the concentration of the population size on the carrying capacity established in \cite{F25}.
Section~\ref{sec:Lookdown} introduces the lookdown construction of both the pre-limiting model and the limiting particle system. In Section~\ref{sec:Conv} we show the forwards in time convergence of the lookdown particle system (relegating some calculations to the Appendix), complete the proof of our main result Theorem \ref{thm:main} and sketch a proof of Theorem \ref{Thm:raphael}. Appendix~\ref{Ap:Lookdown} shows that the proposed lookdown model coincides with the branching process with competition. The stochastic averaging theorem from \cite{F25} is recalled in Appendix~\ref{averagingappendix}. Finally, Appendix~\ref{ratesappendix} contains technical details of the proofs from Section~\ref{sec:Conv}.

\section{Modified lookdown construction}\label{sec:Lookdown}
\begin{subsection}{Modified lookdown construction for the prelimit}
In this section we present a modified lookdown construction for the logistic branching process. Our ordered particle system follows the type-II model of \cite{DK99}, which we now recall. At time \(t\geq 0\), the  state of the population is recorded as a vector
\[
X^{K}(t) = (X_1(t), \ldots, X_{N^{K}(t)}(t)) \in \bigcup_{n \in \mathbf{N} }E^n,
\]
and it evolves according to three types of behaviour. Writing \(N(t)\) for the number of particles in the system at time \(t\). Then:
\begin{itemize}
    \item \textbf{Birth events.} At instantaneous rate \(b p_\ell N\), choose \(\ell+1\) indices uniformly without replacement from \([N+\ell]\), say \(j_1 < \dots < j_{\ell+1}\). The particle with the smallest index, \(j_1\), becomes the parent. New particles with type \(X_{j_1}(t-)\) are inserted at positions \(j_2, \ldots, j_{\ell + 1}\). The rest of the population is reindexed, preserving the original order.  
    
    \item \textbf{Death events.} At instantaneous rate \(N(d+\frac{c}{K}N)\) the particle with the highest index is removed from the population.
    \item \textbf{Mutation.} Each individual mutates according to the generator $\mathcal{Q}^K$, independently from the others.
\end{itemize}
Let \(N = N(x) \coloneqq |x|\) be the current number of particles in the system, The generator of \(X^K(t)\) acts via
\begin{align}
    \label{eq:lookdown_gw}
    \widetilde{\mathcal{L}}_Kf(x) = & \sum_{\ell=1}^\infty\sum_{\substack{J\subset[N+\ell]\\ |J| = \ell + 1}} \frac{bp_\ell N}{\binom{N+\ell}{\ell+1}}
    \left(f(B_{J}x) - f(x)\right) +N\left(d+\frac{c}{K}N\right)\left(f(Dx)) - f(x)\right) 
     \\\notag & +\sum_{i=1}^N\mathcal{Q}_i^Kf(x),
\end{align}
where \(f :\bigcup_{n\in\N}E^n \to \R\) is a bounded function, \(B_{J}\) is the operator that inserts copies of \(x_{j_1}\) at the positions given by the indexes \(J\setminus \{j_1\}\), and re-indexing the remaining elements of \([N + \ell] \setminus J\) while preserving their original order, \(D(x_1, \ldots, x_n) = (x_1, \ldots, x_{n-1})\) is the operator that removes the last particle from the list and \(\mathcal{Q}_i^K\) is the mutation operator \(\mathcal{Q}^K\) acting on the \(i\)-th coordinate of \(f\). Note that \(\lvert B_Jx\rvert = \lvert x \rvert + \lvert J \rvert\).

\begin{figure}
    \centering
    \begin{subfigure}{0.4\textwidth}
        \centering
        \includegraphics[width=\linewidth]{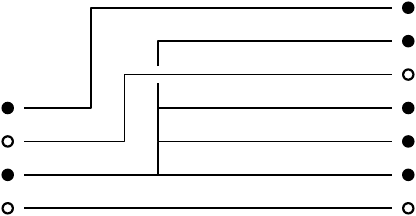}
        \caption{Birth event}
        \label{fig:sub1}
    \end{subfigure}
    \hfill
    \begin{subfigure}{0.4\textwidth}
        \centering
        \includegraphics[width=\linewidth]{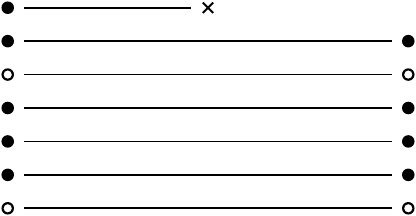}
        \caption{Death event}
        \label{fig:sub2}
    \end{subfigure}
    \caption{In \ref{fig:sub1}, individuals born during the birth event from the second individual are placed in the second, third, forth and fifth level, the rest of individuals are pushed upwards retaining their original order. In \ref{fig:sub2} the individual in the top level is removed from the population. }
    \label{fig:LD}
\end{figure}
\begin{theorem}
    \label{Lookdown} 
    The process of empirical measures \(\left(\sum_{i = 1}^{N^K(t)}\delta_{X^K_i(t)}, t \geq 0\right)\) has the law of the logistic branching process \(\nu^K\) of Definition \ref{DEF:logisticbranchingprocess}. Moreover, if the random vector \(\left(X^K_1(0),\ldots,X^K_{N^K(0)}(0)\right)\) has exchangeable law, then for all \(t>0\), the random vector \(\left(X^K_1(t),\ldots,X^K_{N^K(t)}(t)\right)\) has exchangeable law.
\end{theorem}
This idea is the core of \cite{DK99}. We will show a `modern' proof via the Markov Mapping Theorem \cite[Theorem A.2]{E&K} developed by the same author, although it can be done more directly; we refer the reader to Appendix~\ref{Ap:Lookdown}. This construction of the logistic branching process may seem unusual at first. However, parent–offspring relationships can be inferred directly from the times and subsets \(J\) involved in lookdown birth events. This allows us to select a random sample from the population and trace its genealogy by examining the backward evolution of the model. Moreover, it suffices to characterize the genealogical process obtained by observing only the evolution of the particles at the lowest \(k\) levels. This crucial observation is formalized in the following lemma.
\begin{lemma}\label{sampling}
    Conditional on the event \(N^K(r_KT) \geq k\), the process \(\Pi^K\) has the same law as the coalescent induced by the individuals at the lowest \(k\) levels in the lookdown construction described above. 
\end{lemma}
\begin{proof}
    This follows from the exchangeability of the types established in Theorem \ref{Lookdown}. See \cite[Section 3.10]{E&K} and \cite[Section 5]{DK99} for more details.
\end{proof}
\end{subsection}
\subsection{\texorpdfstring{Limiting \(\boldsymbol{E^\infty}\)-valued lookdown process}{Limiting lookdown process}}

In this section we describe the limiting processes \(X\) with values in \(E^{\infty}\) to which the ordered particle systems \(X^K\) from the previous section converge. These limits are the particle systems providing a genealogical construction of the \(\Lambda\)-Fleming–Viot processes \cite{DK99}. We recall the rigorous construction using families of Poisson point processes. For a finite measure \(\Lambda\) on \([0,1]\) make the decomposition \(\Lambda = a \delta_0 + \Lambda_0\), where \(a \coloneqq \Lambda(\{0\})\) and \(\Lambda_0(\{0\}) = 0\). Let \((\mathfrak{N}_{i,j})_{i<j}\) be independent Poisson processes with rate \(a\), and let \(\mathfrak{M}\) be a Poisson point process on \(\mathbb{R}^+ \times [0,1] \times [0,1]^{\mathbb{N}}\) with intensity \(dt \otimes \frac{\Lambda_0(du)}{u^2} \otimes \left(\ind_{[0, 1]}(v)dv\right)^{\otimes \mathbb{N}}\). Initialize the system with \(X(0)=\left(X_1(0),X_2(0),\ldots\right)\) according to an exchangeable distribution on \(E^\infty\). We build the process \(X\) using these ingredients as follows.
\begin{itemize}
    \item \textbf{Small reproduction events.} For each atom \(t\) in \(\mathfrak{N}_{i,j}\), a particle is inserted at level \(j\) with type \(X_i(t-)\), and the genetic types that were previously at level \(j\) and above are pushed one level up. 
    \item \textbf{Large reproduction events.} For each atom \((t, u, (u_i)_{i\geq1})\in \mathfrak{M}\) define \(J = \{i : u_i \leq u\}\). Now insert particles with type \(X_{\min J}(t-)\) at each level in \(J \setminus \{\min J\}\) and re-index the other particles maintaining their original order.
\end{itemize}
In addition, each particle mutates independently according to the generator \(\mathcal{Q}\), regardless of its level.

\begin{definition}The particle system \((X_t)_{t \geq 0}\) described above is the modified lookdown representation of the \((\Lambda, \mathcal{Q})\)-Fleming–Viot process. In the special case in which \(\Lambda = a\delta_0\) this is the modified lookdown representation of the Fleming-Viot process.
\end{definition}
    
\begin{theorem}[Special case of \cite{Xi09} Theorem 1.1]
    \label{Thm:LFV}
    Let \(X\) be the modified lookdown representation of a \( (\Lambda,\mathcal{Q})\)-Fleming-Viot process. If \(X(0)\) is exchangeable, then \(X(t)\) is exchangeable for all \(t\geq 0\). Moreover, a c\`adl\`ag modification of the process \((\rho(t), t \geq0)\) defined by the de Finetti measure \(\rho(t) \defeq \lim_{\ell\to \infty}\frac{1}{\ell}\sum_{i = 1}^{\ell} \delta_{X_i(t)}\) has the law of the \((\Lambda, \mathcal{Q})\)-Fleming-Viot process.
\end{theorem}
\begin{remark} To see that the above construction can be carried out, note that the type \(X_n\) of a particle at level \(n\) can only change when at least \textit{two} individuals at or below level \(n\) are involved in a reproduction event. The rate at which this occurs is \(\binom{n}{2}a + \int_0^1\P\left(\text{Bin}(n, u) \geq 2\right)\frac{\Lambda(du)}{u^2}\). Note that as \(u \to 0\), \(\P\left(\text{Bin}(n, u) \geq 2\right) \sim \binom{n}{2}u^2 + o(u^3)\). This rate is therefore finite since \(\Lambda\) is a finite measure.
\end{remark}
This construction now comes with the benefit that even in the limiting setting there is a precise notion of genealogy. Each particle in this model represents a line of descent. 
By Theorem \ref{Thm:LFV} the genealogy of a uniformly random sample of size \(k\) is the same as that of the particles at the lowest \(k\) levels. From this representation we can immediately read off the \(\Lambda\)-coalescents. For example consider the case \(\Lambda(\{0\}) = 0\). When there is an atom with second coordinate \(u\), each subset of size \(j\) of any collection of \(k\) levels is involved in the birth event with probability \(u^j(1-u)^{k-j}\), yielding the total rate of such events \(\int_0^1u^j(1-u)^{k-j}\frac{\Lambda(du)}{u^2}\). This is precisely the collision rate of the \(\Lambda\)-coalescent.

\begin{section}{Convergence of the lookdown processes}\label{sec:Conv}
In what follows, embed the prelimiting processes \(X^K\) in \(E^\infty\) by extending the vector as per \cite[Section 3.3]{DK99}. We specify a convention for concreteness. Let \[\beta^K_j(t) = \sup\{s< t:N^K(s) \geq j\}.\] When \(j > N^K(t)\), set \begin{align*}X^K_j(t) =\begin{cases}
     X^K_j(0) &\text{if } \max_{s\leq t}N^K(s) < j,
    \\X_j^K(\beta^K_j(t)-) &\text{otherwise.}
\end{cases} 
\end{align*} In words, particles that no longer correspond to living individuals keep the last type they had while they did, or their initial type if they were never part of the living population. In each subsection here impose the corresponding conditions and notation of Theorems \ref{Thm:raphael} and \ref{thm:main}. Throughout, we let \[\mathcal{D}(A) = \left\{f \in C_b(E^\infty) : \exists k, f = f(x_1, \ldots, x_k)\right\}\] be the continuous (and hence bounded) functions depending on only finitely many indices. \(\mathcal{D}(A)\) will appear as the domain of three slightly different operators \(A\) in each of the three following subsections. Specifying an arbitrary element \(f \in \mathcal{D}(A)\) always comes with a choice of \(k = k(f)\), so it should be understood that every appearance of \(k\) in this section depends on the choice of \(f\). Throughout this section \(E^\infty\) is given the product topology. In what follows, we carry out the proof without mentioning the mutation mechanism, since it is irrelevant for the main result and its inclusion adds only notational burden.

\begin{remark}
    We use the same letter \(k\) to refer to the size of the sample in our genealogies and the support of the test functions since the genealogy of a sample of size \(k\) is precisely the information captured by test functions based on the lowest \(k\) levels.
\end{remark}

\begin{subsection}{Finite variance regime: Kingman's coalescent}
Recall \(n^K(t) = \frac{N^K(Kt)}{K}\) and that in this regime we rescale time by \(K\). We write \(Z^K(t)~=~X^K(Kt)\). We need to know that this rescaled population size concentrates on \(n_* = \frac{bm-d}{c}\). For \(\varepsilon > 0 \), set \[\tau_K = \inf\{t \geq 0: \left|n^K(t) - n_*\right| > \eps\}.\]

\begin{lemma}{(\cite{F25}, Theorem 2)} For all \(\eps > 0\) and \(t\geq 0\), \(\lim_{K \to \infty}\P\left(\tau_K \leq  t\right) = 0\). 
\end{lemma}

\begin{proof}[Proof of Theorem \ref{thm:fwd}, regime 1.]
Let \(f\in \mathcal{D}(A)\), so that \(f = f(x_1, \ldots, x_k)\) for some \(k = k(f)\in \N\).
When the rescaled population size is \(n\) and time has been rescaled by \(K\), then for each \(J \subset [k]\)  with \(|J| = j \geq 2\), the total rate at which precisely the levels in \(J\) undergoes a lookdown birth event (with no levels in \([k]\setminus J\) selected is
\begin{align*}
    R_j(K, n) &\defeq K^2nb\sum_{\ell:\ell\geq j-1} p_{\ell} \frac{\binom{Kn+\ell-k}{\ell + 1 - j}}{\binom{Kn+\ell}{\ell + 1}} \\&= K^2nb\sum_{\ell:\ell\geq j-1}p_\ell\frac{(Kn+\ell-k)!(Kn-1)!(\ell+1)!}{(Kn-k+j-1)!(\ell+1-j)!(Kn+\ell)!}.  
\end{align*} 
Indeed, for the fixed subset \(J\) to be selected we must have a birth of size at least \(j-1\). When there is a birth of size \(\ell \geq j-1\) it must be that the set of \(\ell + 1\) levels chosen from \([Kn + \ell]\) intersects \([k]\) precisely at the levels in \(J\). There are \(\binom{Kn+\ell}{\ell+1}\) ways to select the set of \(\ell + 1\) levels in total, and \(\binom{Kn+\ell-k}{\ell+1-j}\) in which the only levels from \([k]\) chosen are those in \(J\).
We suppress the dependence on \(k\) from the notation. Note that for \(K\) large enough we will not observe any death events amongst the first \(k\) level processes (at least up to time \(\tau_K\)) since deaths always remove the individual with the highest level, which is much greater than \(k\). We will shortly show that \(R_j(K, n) \overset{K\to \infty}{\to} R_j(n)\) uniformly on \(n\in [n_*-\eps, n_*+\eps]\), where \begin{align*}
    R_j(n) = \begin{cases}
        \frac{b(m+m^{(2)})}{n} &\text{if j = 2}\\
        0 &\text{otherwise}.
    \end{cases}
\end{align*} Before doing so we first complete the proof. Note that \[f(Z^K(t\wedge\tau_K))-f(Z^K(0))- \int_0^{t\wedge \tau_K}\sum_{J\subset [k]}R_{|J|}(K, n^K(s))\left(f(B_JZ^K(s)) - f(Z^K(s))\right)ds\] is a martingale with respect to \(\mathcal{F}^K_t = \sigma\left(Z^K_i(s), n^K(s) : i \in [k], s \leq t\right)\) by the previous observation. The tightness of the sequences of processes \(\{f(Z^K_{\cdot \wedge \tau_K})_{t\in[0, T]}, K\geq1\}\) for each \(f\) follows from \cite[Theorem 1]{aldous1978stopping} since the number of possible transitions from any state is bounded by \(2^k\), the transition rates are uniformly bounded on \([n_*-\eps, n_*+\eps]\) as \(K \to \infty\), and \(\norm{f}_\infty < \infty\). This in turn implies tightness of the family \(Z^K\) in \(D_{E^\infty}[0, T]\). To conclude from here one must show that every subsequential limit solves the martingale problem for \(X\), the limiting lookdown process described above. This can be done by a standard application of \cite[Theorem 4.8.10]{ethierkurtz}. We omit the details and turn to the convergence of the rates \(R_j(K, n)\) described above. First, for \(j \geq 3\) and \(n \in [n_* - \eps, n_* + \eps]\), we have by the inequalities \begin{align*}
    \frac{(\ell+1)!}{(\ell + 1 -j)!} &\leq (\ell+1)^j 
    \\ \frac{(Kn+\ell)!}{(Kn+\ell-k)!} &\geq (Kn + \ell-k)^k
\end{align*}that \begin{align*}
    R_j(K, n) \leq K^2(n_* + \eps)b\sum_{j-1 \leq \ell \leq \eps K}p_\ell\frac{1}{(Kn+\ell-k)^k}(Kn)^{k-j}(\ell+1)^j + K^2(n_*+\eps)b\sum_{\ell > \eps K}p_\ell
\end{align*} Recall that we are working under the assumption that there exist a \(\delta>0\) such that \(\sum_{\ell\geq 1}p_{\ell}\ell^{2+\delta}<\infty\), without loss of generality we can assume \(\delta < 1\). For the first term, note that on \(n \in [n_*-\eps, n_* + \eps]\) (and allowing the constant \(C\) to change value from line to line) \begin{align*}
    &K^2 (n_*+\eps)b\sum_{j-1 \leq \ell \leq \eps K}p_\ell\left(\frac{Kn}{Kn + \ell - k}\right)^{k-j}\left(\frac{\ell+1}{Kn+\ell-k}\right)^j \\ &\leq CK^2\frac{1}{K^j}\sum_{j-1 \leq \ell \leq \eps K}p_\ell \ell^j  \\&= CK^2\frac{1}{K^j} \sum_{j-1 \leq \ell \leq \eps K} p_\ell \ell^{2+\delta}\ell^{j-2-\delta}\\ &\leq CK^2\frac{1}{K^j} (\eps K)^{j-2-\delta}\sum_{\ell\leq \eps K} p_\ell \ell^{2+\delta} \leq \frac{C}{K^{\delta}}.\end{align*} For the second term, let \(Z\) be a random variable with law \((p_\ell)\). Then \[K^2(n+\eps)b\P\left(Z > \eps K\right) = K^2(n+\eps)b\P\left(Z^{2+\delta} > (\eps K)^{2+\delta}\right)\leq \frac{C}{K^\delta}\E\left[Z^{2+\delta}\right]\] so on this interval and for \(j\geq 3\), \(R_j(K, n) \leq \frac{C}{K^\delta}\) as \(K \to \infty\). 

For the case \(j = 2\), note that the sum defining \(R_j(K, n)\) now starts from 1, and we have
\begin{align*}\left|\frac{b}{n}\sum_{\ell\geq1}p_\ell \ell(\ell+1)-R_2(K, n) \right| \leq \frac{b}{n_*-\eps}\sum_{\ell\geq1}p_\ell \ell(\ell+1)\left|1 - \frac{(Kn  )^2(Kn+\ell-k)!(Kn-1)!}{(Kn+\ell)!(Kn-k+1)!}\right|\end{align*} and this upper bound goes to \(0\) uniformly for \(n \in [n_*-\eps, n_*+\eps]\) by the dominated convergence theorem.
\end{proof}

\end{subsection}

\begin{subsection}{\texorpdfstring{The \(\boldsymbol{\alpha}\)-stable regime: Beta coalescent}{Beta coalescent regime}}
We now address the second regime of Theorem \ref{thm:fwd} for which we make use of the stochastic averaging approach of \cite{avg}, in particular its adaptation to the present setting in \cite{F25} stated here as Theorem \ref{averagingtheorem}. The difficulty is that the model experiences a separation of timescales effect: the timescale on which \(n^K\) equilibrates around \(n_*\) is much shorter than the rate of the genetic drift (i.e. the process of types). The upshot is that in this setting the fluctuations in the population size `average out' to an effective population size, which is what allows Forien to obtain \(\Lambda\)-Fleming-Viot processes from this model. We relegate the relevant definitions and the statement of Theorem \ref{averagingtheorem} to Appendix \ref{averagingappendix}. Again recall in this regime that \[n^K(t) = \frac{N^K(K^{\alpha-1}t)}{K}\] and that we are rescaling time by \(K^{\alpha-1}\). As remarked in \cite{F25} Section 3.2 one may assume, without loss of generality, that there are constants \(c_0 < n_*\) and \(C_0 > 0\) such that \(c_0 \leq n^K(0) \leq C_0\) a.s.. Set \[\tau_K = \inf\{t > 0 : n^K(t) < \frac{c_0}{2}\}.\] To apply Theorem \ref{averagingtheorem} we require some knowledge of the stopping time \(\tau_K\) and the rescaled population size \(n^K\).
\begin{lemma}{(\cite{F25}, Lemma 3.3)}\label{alphan} One has~\(\sup_{K\geq 1}\sup_{t\geq 0}\E\left[n^K(t)\right]< \infty\).
\end{lemma}
\begin{lemma}{(\cite{F25}, Lemma 3.4)}\label{alphadiv}
    We have that for all \(t \geq 0\) that \(\P\left(\tau_K \leq t\right) \to 0\) as \(K \to \infty\).
\end{lemma}
With the above in hand, the bulk of the remaining work in this section is contained in the following lemma. 
\begin{lemma}\label{alpharate}
    Let \(k,j \in \N\) with \(k\geq j \geq 2\), \(\alpha \in (1, 2)\) and \(c_0 > 0\). Define \[R_j(K, n) = K^\alpha nb\sum_{\ell : \ell\geq j-1} p_\ell \frac{\binom{Kn + \ell - k}{\ell + 1 - j}}{\binom{Kn+\ell}{\ell+1}}\] and \[R_j(n) = \frac{p_0b}{n^{\alpha-1}}\int_0^1u^j(1-u)^{k-j}\frac{u^{1-\alpha}(1-u)^{\alpha-1}}{u^2}du.\] Then \(R_j(K, n) \to R_j(n)\) as \(K \to \infty\) uniformly on \(n \geq \frac{c_0}{2}\).
\end{lemma}
\begin{proof}
    This is elementary but a little fiddly and unenlightening, so we relegate it to Appendix \ref{ratesappendix}.
\end{proof}

\begin{proof}[Proof of Theorem \ref{thm:fwd}, regime 2.]
    As explained above, we will be applying Theorem \ref{averagingtheorem}. Again denote \(Z^K(t) = X^K(K^{\alpha-1}t)\). Let \(R_j(K, n)\) and \(R_j(n)\) be as in Lemma \ref{alpharate}. For \(f \in \mathcal{D}(A)\) (which comes with a choice of \(k\) as explained above), define the operator 
    \[Af(Z, n) = \sum_{J \subset [k]}R_{|J|}(n)\left(f(B_J Z) - f(Z)\right).\] Since \(R_{|J|}(n)\) is decreasing in \(n\) one has that \begin{align}\label{alphaabound}\sup_{t \leq T\wedge\tau_K}|Af(Z^K(t), n^K(t))|\leq 2\norm{f}R_{|J|}(c_0/2).\end{align} Note \(\mathcal{D}(A)\) is dense in \(C_b(E^\infty)\). Similarly to the proof of Theorem \ref{thm:fwd} part 2 we have that 
    \begin{align*}
        f\left(Z^K(t\wedge\tau_K)\right) - f\left(Z^K(0)\right) - \int_0^{t\wedge\tau_K}\sum_{J \subset [k]} R_{|J|}(K, n^K(s))\left(f\left(B_{J}Z^K(s)\right) - f\left(Z^K(s)\right)\right)ds
    \end{align*} is a martingale with respect to the filtration \(\mathcal{F}^K_t = \sigma\left(Z_i^K(s), n^K(s) : i \in [k], s \leq t\right)\). Rewrite this martingale as \begin{align}\label{alphaerror}
        f\left(Z^K(t\wedge\tau_K)\right) - f\left(Z^K(0)\right) - \int_0^{t\wedge\tau_K}Af(Z^K(s), n^K(s))ds + \eps^K(t)
    \end{align} where \begin{align*}
        \eps^K(t) = \sum_{J\subset [k]}\int_0^{t\wedge\tau_K}\left(R_{|J|}(K, n^K(s)) - R_{|J|}(n^K(s))\right)\left(f\left(B_{J}Z^K(s)\right) - f\left(Z^K(s)\right)\right)ds.
    \end{align*} We then bound 
    \begin{align*}\sup_{t\in[0, T]}\left|\eps^K(t)\right| &\leq \sum_{J\subset [k]}2T\norm{f}_\infty\sup_{s\in [0, t\wedge\tau_K]}\left|R_{|J|}(K, n^K(s)) - R_{|J|}(n^K(s))\right| \\&\leq 2^{k+1}T\norm{f}_\infty\sup_{n\geq \frac{c_0}{2}, j \leq k}\left|R_{j}(K, n) - R_{j}(n)\right|,\end{align*} which goes to \(0\) as \(K \to \infty\) (deterministically) by Lemma \ref{alpharate}. We are now almost done: define the measure \(\Gamma^K\) on \(\R^+\times \R^+\) by \[\Gamma^K([0, t]\times B) = \int_0^{t\wedge \tau_K}\ind_{\{n^K(s)\in B\}}ds.\] In the notation of Appendix \ref{averagingappendix} observe that \(\Gamma^K\) is in \(\mathcal{M}^{\R^+}\). Condition 1 of Theorem \ref{averagingtheorem} holds trivially since \(E^\infty\) is compact. Condition 2 was established in Equations \eqref{alphaabound} and \eqref{alphaerror}. Condition 3 follows from Lemma \ref{alphan}. Condition 4 was noted above. Condition 5 is Lemma \ref{alphadiv}. Then by Theorem \ref{averagingtheorem} the pair \((Z^K, \Gamma^K)\) is tight and for every limit point \((X, \Gamma)\) there is a filtration \(\{\mathcal{G}_t\}\) such that 
    \[f(X(s)) - f(X(0)) - \int_{[0, t]\times \R^+}\sum_{J\subset[k]}R_{|J|}(n)\left(f(B_{J}X(s)) - f(X(s))\right)\Gamma(ds, dn)\] 
    is a \(\{\mathcal{G}_t\}\)-martingale. Furthermore, as shown at the end of \cite[Section 4.1]{F25}, every limit point \(\Gamma\) of \(\Gamma_K\) has \(\Gamma(ds, dn) = \delta_{n_*}(dn)ds\) almost surely. This yields that any limit \(X\) renders \[f(X(s)) - f(X(0)) - \int_0^t\sum_{J\subset[k]}R_{|J|}(n_*)\left(f(B_{J}X(s)) - f(X(s))\right)ds\] a \(\{\mathcal{G}_t\}\)-martingale. This martingale problem admits the modified lookdown representation of the \(\Lambda\)-Fleming-Viot process as its unique solution. 
\end{proof}
    
\end{subsection}

\begin{subsection}{Neveu’s branching process with logistic competition: Bolthausen-Sznitman coalescent}

The proof of the final case follows the same pattern as the previous, although the offspring distribution does not have finite first moment so sharper estimates are required for the tightness of the population process. These have already been obtained in \cite{F25}, so we recall those results here. In this section set \[n^K(t) = \frac{1}{K\log(K)}N^K(t)\] and \(n_* = \frac{bp_0}{c}\). As before, one may assume without loss of generality that \(c_0 \leq n^K(0) \leq C_0\) a.s., for some \(c < n_*\) and \(C_0 > 0\). Set again \[\tau_K = \inf\left\{t \geq 0 : n^K(t) < \frac{c_0}{2}\right\}.\] For \(\eps < n_*\), introduce \(V_\eps : (0, \infty) \to \R\) defined by \[V_\eps (n) = \frac{n}{n_* + \eps} - 1 - \log\left(\frac{n}{n_* - \eps}\right).\] The following lemmas establish the tightness of the family \(\left\{n^K(t\wedge \tau_K) : t\geq 0, K \geq 1\right\}\) and the divergence of the times \(\tau_K\). These are quite noteworthy since in this regime \(n^K(t)\) no longer has finite first moment. 

\begin{lemma}[\cite{F25} Lemma 4.1]\label{neveutau} For all \(t\geq 0\), \(\lim_{K\to \infty}\P\left(\tau_K \leq t\right) = 0\).
    
\end{lemma}

\begin{lemma}[\cite{F25} Lemma 4.3]\label{neveutight} For all \(T > 0 \) and \(\eps > 0\) such that \(\inf_{n \geq 0}V_\eps(n) > -1\), one has \[\sup_{K\geq 1}\sup_{t \in [0, T]}\E\left[\log\left(1+V_\eps(n^K(t\wedge \tau_K)\right)\right] < \infty.\]
\end{lemma}

Just as before, we verify the convergence of the rates at which specified subsets of levels are involved in lookdown events.
\begin{lemma}\label{bolthausenrate}
    Let \(k,j \in \N\) with \(k\geq j \geq 2\) and \(\frac{c_0}{2} > 0\). Define 
    \begin{align*}
      R_j(K, n) = K\log(K) nb\sum_{\ell : \ell\geq j-1} p_\ell \frac{\binom{K\log(K)n + \ell - k}{\ell + 1 - j}}{\binom{K\log(K)n+\ell}{\ell+1}}\quad\text{and}\quad R_j = \int_{0}^1u^j(1-u)^{k-j}\frac{bp_0du}{u^2}.  
    \end{align*}
    Then \(R_j(K, n) \to R_j\) as \(K \to \infty\) uniformly on \(n \geq \frac{c_0}{2}\).
\end{lemma}
\begin{proof}
    See Appendix \ref{ratesappendix}.
\end{proof}

\begin{proof}[Proof of Theorem \ref{thm:fwd}, regime 3] With Lemmas \ref{neveutau} and \ref{neveutight} in hand, the remainder of the proof follows the pattern of the previous subsection. Define the operator \(A\) in the same way as before, and identify the same error term. This will converge to 0 deterministically by Lemma \ref{bolthausenrate}. Define again the occupation measure \(\Gamma^K([0, t]\times B) = \int_0^{t\wedge\tau_K}\ind_{\{n^K(s) \in B\}}ds\). The main change is checking conditions 3 and 5 in the Stochastic Averaging Theorem \ref{averagingtheorem} to see that the pair \((X^K, \Gamma^K)\) is tight; these follow from Lemmas \ref{neveutau} and \ref{neveutight}. By the discussion at the end of Section 4.1 in \cite{F25}, the unique subsequential limit point of \((\Gamma^K)_{K \geq 1}\) is \(\Gamma(ds, dn) = \delta_{n_*}(dn)ds\). That yields that for all \(f \in \mathcal{D}(A)\), one has that \[f(X(s)) - f(X(0)) - \int_0^t\sum_{J\subset[k]}R_{|J|}(n_*)\left(f(B_{J}X(s)) - f(X(s))\right)ds\] is a martingale in some filtration \(\mathcal{G}\). Note that this martingale problem \((A, \mathcal{D}(A))\) has the lookdown representation \(X\) as its unique solution.
\end{proof}
\end{subsection}
\end{section}

\begin{section}{Deduction of Theorem \ref{Thm:raphael} and Theorem \ref{thm:main}}
We now complete the proof of our main result, Theorem \ref{thm:main}, concerning the genealogy of the logistic branching process. We will also briefly sketch how Theorem \ref{Thm:raphael}, the main result of \cite{F25}, could also be deduced. We emphasise that the proof of Theorem \ref{thm:fwd} relies heavily on the approach and intermediate lemmas of \cite{F25}, but it is of interest to see how the convergence of the (measure-valued) type frequency processes could be proved via the lookdown particle system. 

\begin{proof}[Proof of Theorem \ref{thm:main} and Theorem \ref{thm:spatial}] To establish the convergence of the genealogy we appeal to the convergence of the times and levels involved in the lookdown events, which encodes the genealogical structure. To that end we require a minor adaptation of Theorem \ref{thm:fwd} since the process \(X^K\) does not encode full information about the lookdown events: a particle can be involved in an event but not change its type, so this cannot be detected purely from the sample path \(X^K\). To record all the required information, associate to each subset \(J \subset [k]\) a counting process \(\left(N^K_J(t): t\geq0\right)\) that records the number of times that the set \(J\) of levels undergoes a lookdown event by time \(t\) (doing nothing if the population has less than \(k\) individuals). Let \(r_K\) be either \(K, K^{\alpha-1}\) or \(1\) respectively according to the regime considered, and for each \(t\leq T\) define \[\mathcal{J}^K(t) = \left(N^K_J(r_Kt):J\subset [k]\right) \in \N^{2^k}.\] We have, by a proof identical to that of Theorem \ref{thm:fwd}, that \(\mathcal{J}^K\) converges weakly in \(D_{\N^{2^k}}[0,T]\) to \(\mathcal{N} = \left(N_J:J\subset [k]\right)\), where the \(N_J\) are independent Poisson processes and the rate of \(N_J\) is the rate at which \(\lvert J\rvert\) of \(k\) blocks coalesce in the corresponding \(\Lambda\)-coalescent. Let \(\chi\) be the subset of \(D_{\N^{2^k}[0,T]}\) consisting of paths with no two coordinates increasing simultaneously, and such that only finitely many jumps are recorded across all \(J\). Note that \(\mathcal{N} \in \chi\) a.s.. To each path \(x \in \chi\), which records full ancestral information at all birth events, we denote the associated ancestral partition process as \(\psi(x) \in D_{\mathcal{P}_k}[0,T]\). Since \(\psi:\chi\to D_{\mathcal{P}_k}[0,T]\) is continuous it follows by the Continuous Mapping Theorem that \(\psi(\mathcal{J}^K)\) converges weakly in \(D_{\mathcal{P}_k}[0,T]\) to \(\psi(\mathcal{N})\). From Lemma \ref{sampling}, we obtain \(\left(\Pi^K(r_Kt): t \in [0, T]\right) \overset{(d)}{=} \psi(\mathcal{J}^K)\). Finally, by the sampling consistency of the \(\Lambda\)-coalescent, we deduce that \(\psi(\mathcal{N})\) has the law of the corresponding \(\Lambda\)-coalescent. Theorem \ref{thm:spatial} follows immediately from Theorem \ref{thm:fwd} and \cite[Theorem 1.26]{K25}.
\end{proof}

\begin{proof}[Sketch proof of Theorem \ref{Thm:raphael}]
Let \(r_K\) be \(K, K^{\alpha-1}\) or \(1\) respectively according to the regime of study. Introduce the stopping time \[\sigma_K^\ell = \inf\{t : N^K(r_Kt) \leq \ell\}.\] Then \(\sigma_K^\ell \to \infty\) in probability as \(K\to \infty\). Let \begin{align*}\eta^K(t) = \frac{1}{N^K(r_Kt)}\sum_{i = 1}^{N^K(r_Kt)}\delta_{X_i^K(r_Kt)}\end{align*} be the empirical measure of the time rescaled lookdown particle system at time \(t\) (at least until the extinction time). Crucially, by Theorem \ref{Lookdown}, \[\left(\eta^K(t), t \in [0, T]\right) \overset{(d)}{=} \left(\rho^K(r_Kt), t \in [0, T]\right).\] Therefore to prove the convergence in distribution of \(\rho^K\) it remains to prove that of \(\eta^K\) which, in turn, will follow from that of \(X^K\) as follows. 
Define further approximates \[\eta^{K, \ell}(t) = \frac{1}{\ell}\sum_{i = 1}^\ell \delta_{X^K_i(r_Kt\wedge\sigma_K^\ell)}\] and \[\eta^\ell(t) = \frac{1}{\ell}\sum_{i=1}^\ell\delta_{X_i(t)}.\] For fixed \(\ell\), \(\eta^{K, \ell}\to \eta^\ell\) in distribution as \(K \to \infty\) by Theorem \ref{thm:fwd} since for all test functions \(\phi: E\to \R\) we have \begin{align*}
    \la \phi, \eta^{K, \ell} \ra = \frac{1}{\ell}\sum_{i=1}^{\ell}\phi(X^K_i) \overset{(d)}{\to} \frac{1}{\ell}\sum_{i = 1}^{\ell}\phi(X_i^K)  
\end{align*} by the convergence of \(X^K\) in \(D_{E^\infty}[0, T]\) where \(E^\infty\) has the product topology. 
One must then guarantee uniform approximation of \(\eta^K\) by \(\eta^{K, \ell}\) and of \(\eta\) by \(\eta^{\ell}\) which is essentially the content of \cite[Lemma 3.4 and Lemma 3.5]{DK99}.
\end{proof}
\end{section}

\appendix
\section{Appendix}
\subsection{Modified lookdown}\label{Ap:Lookdown}
\begin{proof}[Proof Lemma \ref{Lookdown}]\label{proof:lookdown}
The core of this lemma is based on the Markov Mapping Theorem see \cite{K&R, E&K}. This theorem gives conditions under which a projection of a Markov process remains a Markov process, and when martingale properties are preserved. In our application we will consider the map \begin{align*}\gamma: \bigcup_{N\in\N}E^N &\to \mathcal{M}_F(E)
\\(x_1, \ldots, x_N) &\mapsto \sum_{i = 1}^N \delta_{x_i}.
\end{align*}In our case, it is enough to show that the generator \eqref{eq:lookdown_gw} of the lookdown particle system coincides with the one of the logistic branching process when averaged over a uniform permutation of the levels. This corresponds to `forgetting the order' of the particle system which is imposed by the lookdown construction, and the Markov Mapping Theorem then implies that \(X^K(t)\) is exchangeable for each \(t\), and that the process \(\gamma(X^K)\) has the law of the logistic branching process \(\nu^K\). To that end let \(\Sigma_N\) be the set of permutations of the first \(N\) integers and define \(\Theta f: \bigcup_{N\in\N}E^N\to \R\) by
    \[\Theta f(x) = \frac{1}{N!}\sum\limits_{\sigma \in \Sigma_{N}} f(x_\sigma), \quad \text{for}\quad x\in E^N\] for \(f\) bounded and \(x_\sigma=(x_{\sigma_1},\ldots,x_{\sigma_N})\). 
    We want to show that the generator of the lookdown process \eqref{eq:lookdown_gw} once averaged coincides with the one of the logistic branching process \eqref{eq:gen_K_BGW}. Since the set of linear combinations of the functions of the form \(\nu \mapsto \exp(-\la \phi, \nu \ra)\) forms a separating subalgebra of \(C(\mathcal{M}_1(E))\) it will be enough to show that 
    \begin{align*}
        \Theta \widetilde{\mathcal{L}}_{K}\left(f(x)\right)  = \mathcal{L}_{K}\left(\Theta f(x)\right)
    \end{align*}
    for functions of the form 
    \begin{align*}
        f(x) = \exp \left( \la \log g, \gamma(x) \ra \right)=\prod_{z\in \gamma(x)}g(z)
    \end{align*} where \(g:E\to [0, 1]\) is continuous. Here we write \(z \in \gamma(x)\) to denote collection of locations of atoms of the measure \(\gamma(x)\) counted with multiplicity. 
    Recall the generators \eqref{eq:gen_K_BGW} and \eqref{eq:lookdown_gw} 
    \begin{align*}
        \notag \mathcal{L}_{K}F_{h,\phi}(\nu)= &  \sum_{\ell=1}^{\infty}\int_{E}bp_\ell\left(h(\la \phi, \nu\ra +\ell\phi(x))-h(\la \phi, \nu\ra)\right) \nu(dx)\\
        & + \left(d+\frac{c}{K}\la \nu,1\ra\right)\int_{E}\left(h(\la \phi, \nu\ra -\phi(x))-h(\la \phi, \nu\ra)\right) \nu(dx)\\
        & + h'\left(\langle \phi, \nu\rangle\right)\langle \mathcal{Q}^K\phi, \nu\rangle\\
        \widetilde{\mathcal{L}}_Kf(x) = & \sum_{\ell=1}^\infty\sum_{\substack{J\subset[N+\ell]\\ |J| = \ell + 1}} \frac{bp_\ell N}{\binom{N+\ell}{\ell+1}}
    \left(f(B_{J}x) - f(x)\right) +N\left(d+\frac{c}{K}N\right)\left(f(Dx)) - f(x)\right)\\
    &+\sum_{i=1}^N\mathcal{Q}_i^Kf(x)
    \end{align*}
where \(F_{h,\phi}(\nu)\coloneqq h(\la \phi, \nu\ra)\). We will compute the average \(\Theta\widetilde{\mathcal{L}}_K\) of the two parts of the generator in turn before adding them together. First notice that 
\begin{align}\label{Theta}
    \Theta f(x) = \frac{1}{N!} \sum_{\sigma\in \Sigma_N} \prod_{j=1}^N g(x_{\sigma_{j}})= \prod_{j=1}^N g(x_{j}) = \exp\left(\la \log g,\gamma(x)\ra\right).
\end{align}
We need to compute \begin{align*}
\Theta&\left(\sum_{\ell=1}^\infty\sum_{\substack{J\subset[N+\ell]\\ |J| = \ell + 1}} \frac{bp_\ell N}{\binom{N+\ell}{\ell+1}}
    \left(f(B_{J}x) - f(x)\right)\right).
\end{align*} Since \(\Theta\) is linear, we begin by noting that 
\begin{align*}
    \Theta\left(\sum_{\substack{J\subset[N+\ell]\\ |J| = \ell + 1}} \frac{bp_\ell N}{\binom{N+\ell}{\ell+1}}
    f(B_{J}x)\right) &=
    \frac{1}{N!}\sum_{\sigma\in \Sigma_N}\frac{bp_\ell N}{\binom{N+\ell}{\ell +1}} \sum_{J\subseteq [N+\ell]}f(B_{J}x_\sigma)\\ & =  \frac{1}{(N-1)!}\frac{bp_\ell}{\binom{N+\ell}{\ell +1}}\sum_{J\subseteq [N+\ell]}\sum_{\sigma\in \Sigma_N} \prod_{j=1}^N g(x_{\sigma_{j}}) \cdot g(x_{\sigma_{j_1}})^\ell \\
    & = \frac{1}{(N-1)!}\frac{bp_\ell}{\binom{N+\ell}{\ell +1}}\prod_{j=1}^N g(x_{j})\sum_{J\subseteq [N+m]} \sum_{\sigma\in \Sigma_N} g(x_{\sigma_{j_1}})^\ell   \\
    & = \frac{1}{(N-1)!}\frac{bp_\ell}{\binom{N+\ell}{\ell +1}}\prod_{j=1}^N g(x_{j}) \sum_{J\subseteq [N+m]} (N-1)! \sum_{i=1}^N g(x_i)^\ell\\
    & = \frac{bp_\ell}{\binom{N+\ell}{\ell +1}}\binom{N+\ell}{\ell+1} \prod_{j=1}^N g(x_{j})\cdot\sum_{i=1}^N g(x_i)^\ell\\
    & = bp_\ell\int_{E} \exp(\la \log g, \gamma(x) \ra+\ell\log g(z) )\gamma(x)(dz).
\end{align*}
Combining this with expression (\ref{Theta}) we obtain 
\begin{align*}
\Theta&\left(\sum_{\ell=1}^\infty\sum_{\substack{J\subset[N+\ell]\\ |J| = \ell + 1}} \frac{bp_\ell N}{\binom{N+\ell}{\ell+1}}
    \left(f(B_{J}x) - f(x)\right)\right) \\ & = \sum_{\ell=1}^{\infty}bp_\ell\int_{E} \exp(\la \log g, \gamma(x) \ra+\ell\log g(z) )-\exp(\la \log g, \gamma(x) \ra )  \gamma(x)(dz)
\end{align*}
In a similar way, noting that 
\begin{align*}
    \frac{1}{N!}\sum_{\sigma\in\Sigma_{N}} f(Dx_{\sigma}) & = \frac{1}{N!}\sum_{\sigma\in\Sigma_{N}} \prod_{j=1}^{N-1} g(x_{\sigma_{j}}) = \frac{(N-1)!}{N!} \sum_{i=1}^N \prod_{j\neq i }g(x_j) \\
    & = \frac{1}{N}\int_{E}\exp(\la \log g,\gamma(x)\ra-\log g(z))\gamma(x)(dz)
\end{align*}
we obtain  
\begin{align*}
    \Theta&\left(N\left(d+\frac{c}{K}N\right)\left(f(Dx)) - f(x)\right)\right)\\ & = \left(d+\frac{c}{K}\la \gamma(x),1\ra\right)\int_{E}\left(\exp(\la \log g,\gamma(x)\ra-\log g(z))-\exp(\la \log g, \gamma(x)\ra)\right)\gamma(x)(dz).
\end{align*}
For the mutations, observe that 
\begin{align*}
    \Theta\left(\sum_{i=1}^N Q_i^Kf(x)\right) & = \frac{1}{N!}\sum_{\sigma \in \Sigma_N} \sum_{i=1}^N Q_i^Kf(x_\sigma)\\
    & = \prod_{j=1}^Ng(x_j) \frac{1}{N!}\sum_{\sigma \in \Sigma_N}\sum_{i=1}^N \frac{Q^Kg(x_{\sigma_i})}{g(x_{\sigma_i})}\\
    & =  \prod_{j=1}^Ng(x_j)\sum_{i=1}^N \frac{Q^Kg(x_{i})}{g(x_{i})}\\
    & =\exp\left(\la \log g,\gamma(x)\ra\right)\la g^{-1}Q^Kg,\gamma(x)\ra .
\end{align*}
By combining the parts of the generator, we conclude that \(\Theta \widetilde{\mathcal{L}}_{K}\left(f(x)\right)  = \mathcal{L}_{K}\left(\Theta f(x)\right)\) as required.

\end{proof}
\subsection{Stochastic averaging}\label{averagingappendix}

For a metric space \(S\), denote by \(\mathcal{M}^S\) the subset of \(\mu \in \mathcal{M}_F(\R^+\times S)\) such that for all \(t\geq0\), \[\mu([0, t]\times S)\leq t.\] Letting \(\mu^t\) denote the restriction of \(\mu \in \mathcal{M}^S\) to \([0, t]\times S\) and \(d^t\) the Prokhorov metric on \(\mathcal{M}_F([0, t]\times S)\), metrize \(\mathcal{M}^S\) by defining \[d(\mu, \nu) = \int_0^\infty e^{-t}(d^t(\mu^t, \nu^t)\wedge 1)dt.\] 
\begin{theorem}[\cite{F25}, Theorem 5. Adapted from \cite{avg}, Theorem 2.1]\label{averagingtheorem}
    Let \(E_1\) and \(E_2\) be two complete separable metric spaces and set \(E := E_1 \times E_2\).
Suppose that for \(K \geq 1\), \(\{(X_K, Y_K), K \geq 1\}\) is a sequence of random variables taking values in \(D([0,\infty) \times E)\), adapted to some filtration \((\mathcal{F}_t^K, t \geq 0)\). We assume the following.

\begin{enumerate}
    \item The sequence \((X_K, K \geq 1)\) satisfies the compact containment condition, i.e. for each \(\varepsilon > 0\) and \(T > 0\), there exists a compact set \(\mathcal{K} \subset E_1\) such that
    \[
    \inf_{K \geq 1} \mathbb{P}\left( X_K(t) \in \mathcal{K}, \, \forall t \in [0,T] \right) \geq 1 - \varepsilon.
    \]
    
    \item There exists an operator \(A : \mathcal{D}(A) \subset C_b(E_1) \to C_b(E_1 \times E_2)\) and a sequence of \(\mathcal{F}_t^K\)-stopping times \((\tau_K, K \geq 1)\) such that, for all \(f \in \mathcal{D}(A)\), there exists a process \((\varepsilon_K^f(t), t \geq 0)\) such that
    \[
    f(X_K(t \wedge \tau_K)) - f(X_K(0)) - \int_0^{t \wedge \tau_K} Af(X_K(s), Y_K(s)) ds + \varepsilon_K^f(t)
    \]
    is an \(\mathcal{F}_t^K\)-martingale, and, for all \(T > 0\), there exists a constant \(C_T > 0\) such that
    \[
    \sup_{t \in [0, T \wedge \tau_K]} |Af(X_K(t), Y_K(t))| \leq C_T, \quad \text{almost surely},
    \]
    and a deterministic sequence \((\eta_K, K \geq 1)\), converging to zero as \(K \to \infty\), such that
    \[
    \sup_{t \in [0,T]} |\varepsilon_K^f(t)| \leq \eta_K.
    \]
    
    \item For any \(T > 0\), the collection of \(E_2\)-valued random variables \(\{Y_K(t \wedge \tau_K), t \in [0,T], K \geq 1\}\) is tight.
    
    \item The set \(\mathcal{D}(A)\) is dense in \(C_b(E_1)\) for the topology of uniform convergence on compact sets.
    
    \item For any \(t \geq 0\),
    \[
    \lim_{K \to \infty} \mathbb{P}(\tau_K \leq t) = 0.
    \]
\end{enumerate}

We then define a sequence \((\Gamma_K, K \geq 1)\) of \(\ell(E_2)\)-valued random variables as
\[
\Gamma_K([0,t] \times B) := \int_0^{t \wedge \tau_K} \ind_{\{Y_K(s) \in B\}} \, ds,
\]
for all measurable \(B \subset E_2\). Then, under these assumptions, the sequence of \(D([0,\infty), E_1) \times \ell(E_2)\)-valued random variables \(\{(X_K, \Gamma_K), K \geq 1\}\) is tight and, for any limit point \((X, \Gamma)\),
\[
\Gamma([0,t] \times E_2) = t, \quad \forall t \geq 0,
\]
almost surely and there exists a filtration \((\mathcal{G}_t, t \geq 0)\) such that, for all \(f \in \mathcal{D}(A)\),
\[
f(X(t)) - f(X(0)) - \int_{[0,t] \times E_2} Af(X(s), y) \Gamma(ds, dy)
\]
is a \(\mathcal{G}_t\)-martingale.

\end{theorem}

\subsection{Convergence of the transition rates}\label{ratesappendix}

\begin{proof}[Proof Lemma \ref{alpharate}]
Make the decomposition 
\begin{equation}
    \notag R_{j}(K, n) - R_j(n) \coloneqq E_1(K,n)+E_2(K,n)+E_3(K,n)
\end{equation}
where 
\begin{align}
    E_1(K,n)  \coloneqq&  K^{\alpha}nb\sum_{\ell\geq j-1}p_\ell\frac{(Kn+\ell-k)!(Kn-1)!(\ell+1)!}{(Kn-k+j-1)!(\ell+1-j)!(Kn+\ell)!}\label{nofactorials} \\
    & \notag  \quad-  K^\alpha nb\sum_{\ell\geq 1}p_\ell\left(\frac{Kn}{Kn+\ell}\right)^{k-j}\left(\frac{\ell}{Kn+\ell}\right)^j
\end{align}
\begin{align}
    E_2(K,n) & \coloneqq K^\alpha nb\sum_{\ell\geq 1}p_\ell\left(\frac{Kn}{Kn+\ell}\right)^{k-j}\left(\frac{\ell}{Kn+\ell}\right)^j\label{sumtoint}\\
    \notag &\quad - nbp_0\int_{\frac{1}{K}}^\infty\frac{1}{y^{1+\alpha}}\left(\frac{n}{n+y}\right)^{k-j}\left(\frac{y}{n+y}\right)^jdy     
\end{align}

\begin{align}
    \label{initialsegment} E_3(K,n) & \coloneqq nbp_0\int_{\frac{1}{K}}^\infty\frac{1}{y^{1+\alpha}}\left(\frac{n}{n+y}\right)^{k-j}\left(\frac{y}{n+y}\right)^jdy \\
    \notag & \quad - nbp_0\int_0^\infty\frac{1}{y^{1+\alpha}}\left(\frac{n}{n+y}\right)^{k-j}\left(\frac{y}{n+y}\right)^jdy  
\end{align}    
To see this, note that by substitution \[nbp_0\int_0^\infty\frac{1}{y^{1+\alpha}}\left(\frac{n}{n+y}\right)^{k-j}\left(\frac{y}{n+y}\right)^jdy = \frac{p_0b}{n^{\alpha-1}}\int_0^1u^j(1-u)^{k-j}\frac{(1-u)^{\alpha-1}u^{1-\alpha}}{u^2} = R_j(n).\] 

We begin with the first term, \(E_1(K, n)\). Notice that it is controlled by standard bounds on factorials, together with the observation that the first \(j-2\) terms of the second sum clearly go to 0 as \(K \to \infty\) since \(\alpha < 2 \leq j\).

For the second term \(E_2(K, n)\), first define \(p(y) = y^{1+\alpha}p_{\lfloor y \rfloor}\), so that \(p(y) \to p_0\) as \(y \to \infty\). Define also \(y_K = \frac{\lfloor Ky \rfloor}{K}\). Then \[K^\alpha nb\sum_{\ell\geq1}p_\ell\left(\frac{Kn}{Kn+\ell}\right)^{k-j}\left(\frac{\ell}{Kn+\ell}\right)^j = \int_{\frac{1}{K}}^\infty \frac{p(Ky)}{y^{1+\alpha}}\left(\frac{n}{n + y_K}\right)^{k-j}\left(\frac{y_K}{n + y_K}\right)^jdy.\] We can then further decompose term \eqref{sumtoint} as 
        \begin{align}&nb\int_{\frac{1}{K}}^\infty \frac{p(Ky) - p_0}{y^{1+\alpha}}\left(\frac{n}{n + y_K}\right)^{k-j}\left(\frac{y_K}{n + y_K}\right)^jdy\label{pK} \\&+ nb\int_{\frac{1}{K}}^1 \frac{p_0}{y^{1+\alpha}}\left(\left(\frac{n}{n + y_K}\right)^{k-j}\left(\frac{y_K}{n + y_K}\right)^j - \left(\frac{n}{n + y}\right)^{k-j}\left(\frac{y}{n + y}\right)^j\right)dy \label{1K1}\\&+ nb\int_1^\infty \frac{p_0}{y^{1+\alpha}}\left(\left(\frac{n}{n + y_K}\right)^{k-j}\left(\frac{y_K}{n + y_K}\right)^j - \left(\frac{n}{n + y}\right)^{k-j}\left(\frac{y}{n + y}\right)^j\right)dy.\label{tail}\end{align} 
        Using the fact that \(j\geq 2\), so that \(\left(\frac{y_K}{n+y_K}\right)^j\leq \left(\frac{y}{n+y}\right)^j \leq \frac{y^2}{(n + y)^2}\leq \frac{y^2}{n^2}\) we can bound the modulus of expression \eqref{pK} above by \(\frac{b}{n}\int_{\frac{1}{K}}^1 \left|\frac{p(Ky)-p_0}{y^{\alpha-1}}\right|dy + b\int_1^\infty \left|\frac{p(Ky)-p_0}{y^\alpha}\right|dy,\) which converges to \(0\) uniformly on \(n\geq \frac{c_0}{2}\) by the dominated convergence theorem.
        As for expression \eqref{1K1}, first rewrite it as \begin{align}&nb\int_{\frac{1}{K}}^1 \frac{p_0}{y^{1+\alpha}} \left(\left(\frac{n}{n + y_K}\right)^{k-j} - \left(\frac{n}{n + y}\right)^{k-j}\right)\left(\frac{y_K}{n+y_K}\right)^j dy\label{bruh} \\&+ nb\int_{\frac{1}{K}}^1 \frac{p_0}{y^{1+\alpha}}\left(\frac{n}{n+y}\right)^{k-j}\left(\left(\frac{y_K}{n+y_K}\right)^j - \left(\frac{y}{n+y}\right)^j\right)dy.\label{ohboy}\end{align} For term (\ref{bruh}), if \(k = j\) it is \(0\). Otherwise its modulus can be bounded above by \[nb\int_{\frac{1}{K}}^1 \frac{p_0}{y^{1+\alpha}}(k-j)\left|\frac{n}{n+y_K}-\frac{n}{n+y}\right|\frac{y^2}{n^2}dy \leq \frac{b(k-j)}{n^2K}\int_{\frac{1}{K}}^1\frac{p_0}{y^{\alpha-1}}dy\] where for the last inequality we use that \[\left|\frac{n}{n+y_K}-\frac{n}{n+y}\right| \leq \frac{|y_K - y|}{n}\leq \frac{1}{nK}.\] This bound converges to \(0\) uniformly on \(n\geq \frac{c_0}{2}\). Term (\ref{ohboy}) requires more care. Since \(j\geq2\), noting for \(a, b <1\) that \(|a^j - b^j| = |(a^2)^{\frac{j}{2}} - (b^2)^{\frac{j}{2}}| \leq \frac{j}{2}|a^2 - b^2|\), we have
        \begin{align*}
            \left|\left(\frac{y_K}{n+y_K}\right)^j - \left(\frac{y}{n+y}\right)^j\right| & \leq \frac{j}{2}\left|\left(\frac{y_K}{n+y_K}\right)^2 - \left(\frac{y}{n+y}\right)^2\right|\\ &\leq \frac{j}{2}\frac{n^2\left|y_K^2 - y^2\right| + 2nyy_K\left|y_K-y\right|}{n^4}.
        \end{align*}
        Then use that \(|y^2-y_K^2|\leq \frac{4y}{K}\) for \(y\geq \frac{1}{K}\) and that \(\lim_{K\to\infty}\frac{1}{K}\int_{\frac{1}{K}}^1 \frac{1}{y^\alpha}dy = 0\), together with the observation \(2nyy_K|y_K-y|\leq \frac{2ny^2}{K}\) to obtain that term (\ref{ohboy}) is bounded above by \begin{align*}
            &nb\int_{\frac{1}{K}}^1\frac{p_0}{y^{1+\alpha}}\frac{j}{2}\frac{\frac{4ny}{K} + \frac{2ny^2}{K}}{n^4}dy \leq \frac{bj}{n^2}\frac{2}{K}\int_{\frac{1}{K}}^1\frac{p_0}{y^\alpha}dy + \frac{b}{n^2K}\int_{\frac{1}{K}}^1\frac{p_0}{y^{\alpha-1}}dy,
        \end{align*} and both terms in this sum converge uniformly in \(n\geq \frac{c_0}{2}\) as \(K \to \infty\).
        
        For a bound on expression \eqref{tail}, write \begin{align}&nb\int_1^\infty \frac{p_0}{y^{1+\alpha}} \left(\left(\frac{n}{n + y_K}\right)^{k-j} - \left(\frac{n}{n + y}\right)^{k-j}\right)\left(\frac{y_K}{n+y_K}\right)^j dy \label{whatislove}\\&+ nb\int_1^\infty \frac{p_0}{y^{1+\alpha}}\left(\frac{n}{n+y}\right)^{k-j}\left(\left(\frac{y_K}{n+y_K}\right)^j - \left(\frac{y}{n+y}\right)^j\right)dy.\label{babydonthurtme}\end{align} For term (\ref{whatislove}), if \(k=j\) this is \(0\). For \(k-j\geq1\) we can bound the modulus of this by \(nb\int_1^\infty \frac{p_0}{y^{1+\alpha}}(k-j) \left|\frac{n}{n+y_K} - \frac{n}{n+y}\right|\left(\frac{y_K}{n+y_K}\right)^j \leq \frac{b(k-j)}{Kn}\int_1^\infty \frac{p_0}{y^{1+\alpha}}\) which converges to \(0\) uniformly on \(n\geq \frac{c_0}{2}\). For term (\ref{babydonthurtme}), since \(j\geq2\) we can bound by \[nb\int_1^\infty \frac{p_0}{y^{1+\alpha}}j\left|\frac{y_k}{n+y_K} - \frac{y}{n+y}\right|dy \leq \frac{bj}{nK}\int_1^\infty \frac{p_0}{y^{1+\alpha}},\] which also converges uniformly.

Finally, for the third term \(E_3(K, n)\), we have \[\left|-nbp_0\int_0^{\frac{1}{K}}\frac{1}{y^{1+\alpha}}\left(\frac{n}{n+y}\right)^{k-j}\left(\frac{y}{n+y}\right)^jdy\right| \leq \frac{bp_0}{n}\int_0^{\frac{1}{K}}\frac{1}{y^{1+\alpha}}y^2dy = \frac{bp_0}{n(2-\alpha)K^{2-\alpha}}\] which goes to \(0\) uniformly on \(n \geq \frac{c_0}{2}\). This handles each of the errors in the original decomposition.
\end{proof}

\begin{proof}[Proof Lemma \ref{bolthausenrate}]
Let \(a_K = K\log(K)\) and define again \(p(y) = y^2p_{\lfloor y\rfloor}\) so that \(p(y)\to \infty\) as \(y \to \infty\). In this section we denote \(y_K = \frac{\lfloor a_Ky\rfloor}{a_K}\). 
\begin{equation*}
    R_{j}(K, n) - R_j(n) = E_1(K,n)+E_2(K,n)+E_3(K,n)
\end{equation*}
\begin{align*}
    E_1(K,n) = & a_Knb\sum_{\ell\geq j-1}p_\ell\frac{(a_Kn+\ell-k)!(a_Kn-1)!(\ell+1)!}{(a_Kn-k+j-1)!(\ell+1-j)!(a_Kn+\ell)!}\\
    \notag &\quad -  a_K nb\sum_{\ell\geq1}p_\ell\left(\frac{n}{n+\ell}\right)^{k-j}\left(\frac{\ell}{a_Kn+\ell}\right)^j, 
\end{align*}
\begin{align*}
    E_2(K,n) = & a_K n b\sum_{\ell\geq 1}p_\ell\left(\frac{a_Kn}{a_Kn+\ell}\right)^{k-j}\left(\frac{\ell}{a_Kn+\ell}\right)^j\\
    \notag& - nbp_0\int_{\frac{1}{a_K}}^\infty \frac{1}{y^2}\left(\frac{n}{n+y}\right)^{k-j}\left(\frac{y}{n+y}\right)^jdy, 
\end{align*} and 
\begin{align*}
    E_3(K,n) & = nbp_0\int_{\frac{1}{a_K}}^\infty\frac{1}{y^{2}}\left(\frac{n}{n+y}\right)^{k-j}\left(\frac{y}{n+y}\right)^jdy \\
    \notag & - nbp_0\int_{0}^\infty \frac{1}{y^2}\left(\frac{n}{n+y}\right)^{k-j}\left(\frac{y}{n+y}\right)^jdy.
\end{align*} Note again that \[nbp_0\int_0^\infty\frac{1}{y^2}\left(\frac{n}{n+y}\right)^{k-j}\left(\frac{y}{n+y}\right)^jdy = \int_0^1(1-u)^{k-j}u^j\frac{bp_0}{u^2}du.\]

For the first term \(E_1(K, n)\).The analysis is the same as before. Similarly for the second term, \(E_2(K, n)\). Note that \[a_K nb\sum_{\ell\geq1}p_\ell\left(\frac{a_Kn}{a_Kn+\ell}\right)^{k-j}\left(\frac{\ell}{a_Kn+\ell}\right)^j = nb\int_{\frac{1}{a_K}}^\infty \frac{p(a_Ky)}{y^2}\left(\frac{n}{n+y_K}\right)^{k-j}\left(\frac{y_K}{n+y_K}\right)^jdy,\] so write this term as \begin{align}
        &nb\int_{\frac{1}{a_K}}^\infty \frac{p(a_Ky) - p_0}{y^2}\left(\frac{n}{n+y_K}\right)^{k-j}\left(\frac{y_K}{n+y_K}\right)^jdy \label{x}
        \\&+nb\int_{\frac{1}{a_K}}^\infty \frac{p_0}{y^2}\left(\left(\frac{n}{n+y_K}\right)^{k-j} - \left(\frac{n}{n+y}\right)^{k-j}\right)\left(\frac{y_K}{n+y_K}\right)^jdy \label{y}
        \\&+nb\int_{\frac{1}{a_K}}^\infty \frac{p_0}{y^2}\left(\frac{n}{n+y}\right)^{k-j}\left(\left(\frac{y_K}{n + y_K}\right)^j - \left(\frac{y}{n + y}\right)^j\right)dy.\label{z}
    \end{align}
We analyse each term individually. Lets start with expression \eqref{x}, which is bounded by \begin{align*}&\frac{b}{n}\int_{\frac{1}{a_K}}^n\left|p(a_Ky)-p_0\right|dy + nb\int_n^\infty\frac{\left|p(a_Ky) - p_0\right|}{y^2}dy
        \\& = b\int_{\frac{1}{na_K}}^1\left|p(a_Kny) - p_0\right|dy + b \int_1^\infty \frac{\left|p(na_Ky)-p_0\right|}{y^2}dy, 
        \end{align*}and both of these terms go to \(0\) uniformly as \(K \to \infty\) by the dominated convergence theorem. For term \eqref{y}, when \(k-j\geq 1\) then we divide \eqref{y} in to two integrals, one from \(\frac{1}{a_K}\) to \(1\) and a second from \(1\) to \(\infty\). The modulus of the first integral can be bounded above by 
        \begin{align*}
            nb\int_{\frac{1}{a_K}}^1 \frac{p_0}{y^2}\left|\left(\frac{n}{n+y_K}\right)^{k-j} - \left(\frac{n}{n+y}\right)^{k-j}\right|&\left(\frac{y_K}{n+y_K}\right)^jdy \\
            & \leq nb(k-j)\int_{\frac{1}{a_K}}^1 \frac{p_0}{y^2}\left| \frac{n}{n+y_K}-\frac{n}{n+y} \right|\frac{y^2}{n^2} dy \\
            & \leq \frac{bp_0(k-j)}{n^2}\int_{\frac{1}{a_K}}^1 \left|y-y_{K}\right| dy\leq \frac{bp_0(k-j)}{n^2a_K}
        \end{align*} where in the first inequality we use \(j\geq 2\). This bound goes to \(0\) uniformly on \(n\) as \(K\) goes to infinity. For the integral from \(1\) to \(\infty\) notice that
        \begin{align*}
            nb\int_{1}^\infty \frac{p_0}{y^2}\left|\left(\frac{n}{n+y_K}\right)^{k-j} - \left(\frac{n}{n+y}\right)^{k-j}\right|\left(\frac{y_K}{n+y_K}\right)^jdy &\leq nb(k-j)\int_{1}^\infty \frac{p_0}{y^2} \frac{1}{na_{K}} dy 
        \end{align*}
        which also converges to \(0\) uniformly on \(n\geq \frac{c_0}{2}\). For the equation (\ref{z}) we again split the integral at \(1\). First note that since \(j \geq 2\)
        \begin{align*}
            &\left|nb\int_{\frac{1}{a_K}}^1 \frac{p_0}{y^2}\left(\frac{n}{n+y}\right)^{k-j}\left(\left(\frac{y_K}{n + y_K}\right)^j - \left(\frac{y}{n + y}\right)^j\right)dy\right| \\
            &\quad \quad\leq nb\int_{\frac{1}{a_K}}^1 \frac{p_0}{y^2}\left|\left(\frac{y_K}{n + y_K}\right)^j - \left(\frac{y}{n + y}\right)^j\right|dy 
            \\ &\quad \quad \leq \frac{nbj}{2}\int_{\frac{1}{a_K}}^1\frac{p_0}{y^2}\left|\left(\frac{y_K}{n + y_K}\right)^2 - \left(\frac{y}{n + y}\right)^2\right|dy \\
            &\quad \quad \leq  \frac{nbj}{2}\int_{\frac{1}{a_K}}^1\frac{p_0}{y^2}\frac{\left|n^2(y_K^2 - y^2)\right| + y_ky\left|2n(y_K-y)\right|}{n^4}dy \\ &\quad \quad\leq \frac{2bj}{na_K}\int_{\frac{1}{a_K}}^1\frac{p_0}{y}dy + \frac{bj}{n^2a_K}p_0.
        \end{align*} Where in the final step we use \(|y_K^2-y^2| \leq \frac{4y}{a_K}\). Note also that\(\frac{1}{a_K}\int_{\frac{1}{a_K}}^1\frac{p_0}{y^2} \to 0\) as \(K \to \infty\).
        Also we have 
        \begin{align*}
            &\left|nb\int_1^\infty \frac{p_0}{y^2}\left(\frac{n}{n+y}\right)^{k-j}\left(\left(\frac{y_K}{n + y_K}\right)^j - \left(\frac{y}{n + y}\right)^j\right)dy\right| \\
            & \quad \quad\leq nb\int_1^\infty \frac{p_0}{y^2}\left|\left(\frac{y_K}{n + y_K}\right)^j - \left(\frac{y}{n + y}\right)^j\right|dy 
            \\& \quad \quad\leq nbj\int_1^\infty\frac{p_0}{y^2}\left|\frac{y_K}{n+y_K}-\frac{y}{n+y}\right|dy \leq \frac{bj}{a_{K}}\int_1^\infty \frac{p_0}{y^2}dy,
        \end{align*} which also goes to \(0\) uniformly on \(n\geq\frac{c_0}{2}\).

    Finally for the third term, \(E_3(K, n)\). We have that 
    \begin{align*}
    \left|-nbp_0\int_0^{\frac{1}{a_K}}\frac{1}{y^{2}}\left(\frac{n}{n+y}\right)^{k-j}\left(\frac{y}{n+y}\right)^jdy\right| \leq \frac{bp_0}{n}\int_0^{\frac{1}{a_{K}}}\frac{1}{y^{2}}y^2dy=\frac{bp_0n}{a_K}            
    \end{align*}
    which goes to \(0\) uniformly on \(n \geq \frac{c_0}{2}\).
\end{proof}
\section*{Acknowledgement}
We are grateful to Alison Etheridge for many lively and instructive discussions on the subject of this article, and to Peter Koepernik for pointing out the convergence to the spatial coalescent and modifying his preprint to incorporate the modified lookdown construction used here. We also thank Rapha\"el Forien and F\'elix Foutel-Rodier for helpful discussions and feedback.

\printbibliography
\end{document}